\documentclass[12pt]{article}

\usepackage{amsfonts,psfig}

\newcounter{main}
\newenvironment{main}{\noindent\refstepcounter{main} {\bf Theorem \Alph{main}: }\em}{}

\newtheorem{theorem}{Theorem}[section]
\newtheorem{proposition}[theorem]{Proposition}
\newtheorem{lemma}[theorem]{Lemma}
\newtheorem{corollary}[theorem]{Corollary}

\newtheorem{remark}{Remark}[section]
\newtheorem{definition}{Definition}[section]
\newtheorem{maintheorem}{Theorem}

\newcommand{\blanksquare}{\,\,\,$\sqcup\!\!\!\!\sqcap$}
\newenvironment{proof}{{\flushleft {\bf Proof: }}}{\blanksquare}
\newenvironment{rproof}[1]{{\flushleft{\bf Proof of #1: }}}{\blanksquare}
\newcounter{example}
\newenvironment{example}%
{{\stepcounter{example}}{\flushleft {\bf Example \arabic{example}:}}}%
{\par}

\newcommand{\al} {\alpha}       
\newcommand{\be} {\beta}        
    
\newcommand{\de} {\delta}       \newcommand{\De}{\Delta}

\newcommand{\ep} {\epsilon}

\newcommand{\la} {\lambda}      \newcommand{\La}{\Lambda}

\newcommand{\si} {\sigma}

\newcommand{\vfi}{\varphi}

\newcommand{\diame}[1]{{\rm diam\,}(#1)}

\newcommand{\dimns}{{\rm dim\,}}
\newcommand{\dista}[2]{{\rm dist\,}(#1,#2)}

\newcommand{\inter}[1]{{\rm int\,}(#1)}

\newcommand{\refer}[1]{(\ref{#1})}

\newcommand{\supp}{{\rm supp\,}}

\newcommand{\natur}{{\mathbb N}}
\newcommand{\real} {{\mathbb R}}
\newcommand{\difsim}{{\bigtriangleup}}

\newcommand{\relativ}{{\mathbb Z}}

\newcommand{\un}{\underline}

\newcommand{\ov}{\overline}
\newcommand{\pa}{\partial}
\newcommand{\til}{\tilde}

\newcommand{\mapto}{\longrightarrow}
\newcommand{\impl}{\Longrightarrow}

\newcommand{\SF}{{\cal F}}

\newcommand{\SM}{{\cal M}}

\newcommand{\SV}{{\cal V}}

\newcommand{\qand}{\quad\mbox{and}\quad}
\newcommand{\subl}{\underline}
\newcommand{\tra}{\overline}

\newcounter{contador}
\newenvironment{nlista}[1]%
{\begin{list}{#1}{\usecounter{contador}}}%
{\end{list}}

{\begin{list}{#1}{\usecounter{#2}}}%
{\end{list}}

{\begin{list}{#1}%
{\usecounter{contador}\setlength{\leftmargin}%
{0.5cm}\setlength{\labelsep}{0.1cm}\setlength{\itemsep}%
{0cm}\setlength{\parsep}{0.1cm}}}%
{\end{list}}

\newcommand{\Letra}{{\underline{\Alph{contador}}.}}

\title{Infinitely Many Stochastically Stable Attractors}

\author{V\'{\i}tor Ara\'ujo\thanks{Partially supported by FCT through
    Centro de Matem\'atica da Universiade do Porto}}

\date{ }

\begin{document}
\maketitle


\begin{abstract}

Let $f$ be a diffeomorphism of a compact finite dimensional
boundaryless manifold $M$ exhibiting infinitely many coexisting
attractors. Assume that each attractor supports a stochastically
stable probability measure and
that the union of the basins of attraction of each attractor
covers Lebesgue almost all
points of $M$. We prove that the time averages of almost all orbits
under random perturbations are given by a finite number of probability
measures. Moreover these probability measures are close to the
probability measures supported by the
attractors when the perturbations are close to the original map $f$.

\end{abstract}

\section{Introduction}
\label{intro}

Let $f:M\mapto M$ be a diffeomorphism of a
compact finite dimensional boundaryless manifold
$M$ exhibiting infinitely many coexisting attractors
--- e.g. the examples provided
first by S. Newhouse
in~\cite{N1,N2,N3},
or later by Pumari\~no-Rodriguez
\cite{Pu} and E. Colli~\cite{Co},
and more recently by Bonatti-D\'{\i}az~\cite{BD}.
What can be said about the stability of time
averages of continuous functions $\varphi:M\mapto\real$
\begin{equation}\label{averagenonoise}
\lim_{n\to\infty} \frac1n\sum_{j=0}^{n-1} \varphi( f^j (z) ),
\qquad z\in M
\end{equation}
under
random perturbations introduced at each iteration?

For general uniformly hyperbolic
diffeomorphisms (Axiom A case)
Sinai~\cite{Sn}, Ruelle~\cite{R} and Bowen~\cite{BR} have shown
that there are finitely many probability measures
$\mu_1,\ldots,\mu_N$ such that
\begin{equation}\label{limitaveragenonoise}
\lim_{n\to\infty} \frac1n\sum_{j=0}^{n-1} \varphi( f^j (z) )
= \int \varphi \,d\mu_i,
\end{equation}
for all continuous $\varphi:M\to\real$ and for points $z$ in a
positive Lebesgue measure subset $B(\mu_i)$ of $M$, $i=1,\ldots, l$.
These measures are called $SRB$ measures and the sets $B(\mu_i)$ are
the \emph{basins} of these measures, $i=1,\ldots, l$.
Moreover the union of these
basins covers Lebesgue almost all points of $M$ and the supports
$\La_i=\supp{\mu_i}$ are attractors, i.e., $\La_i$ is a compact
$f$-invariant transitive set whose basin of attraction is a
neighborhood of $\La_i$.

Recently Bonatti-Viana~\cite{BV2}
together with Alves~\cite{ABV} prove
results analogous to the results of Sinai,
Ruelle and Bowen for some maps with nonuniform
expansion: there are finitely many SRB measures
whose basins cover Lebesgue almost every
point of the ambient manifold.
It is worth noting that the mere \emph{existence}
of time averages for general systems is still
an open problem.

The problem of stability of the measures $\mu_i$ has been studied for
uniformly hyperbolic systems by Kifer~\cite{Ki,Ki2} and
Young~\cite{Y}. It has been proved that if we fix a $\La_i$ and take
at each iteration a diffeomorphism $f_j$ close to $f$,
$j=0,1,2\ldots$, then the time averages (subject to a given
probability distribution for the choices of $f_j$)
\begin{equation}\label{averagewithnoise}
\lim_{n\to + \infty} \frac1n \sum_{j=0}^{n-1}
\varphi\left( f_j\circ\cdots\circ f_1(z) \right),
\qquad z\in M,
\end{equation}
exist for a positive Lebesgue measure subset of points $z\in M$ close
to $\La_i$ and are equal to $\int \varphi \,d\mu_i'$ for a probability
measure $\mu_i'$.
Moreover $\mu_i'$ gets close to the probability measure $\mu_i$ when
the $f_j$ are restricted to a small neighborhood of $f$. This holds
for all $i=1,\ldots, l$.
We say that such hyperbolic systems are \emph{stochastically stable}.
In~\cite{BV} Benedicks and Viana show that the H\'enon map, for a
positive Lebesgue set of parameters, also is stochastically stable.

A probability measure $\mu'$ on $M$ is \emph{physical} if the time
averages~\refer{averagewithnoise} of every continuous
$\varphi:M\to\real$, for almost all
choices of $f_j$, coincide with the space average $\int \varphi
\,d\mu'$ for a positive Lebesgue measure subset of points $z\in M$.

In~\cite{Ar} the author proved that,
under physical random perturbations along
finite dimensional parameterized families of
diffeomorphisms, there exist time averages
for almost every perturbed orbit of each point of $M$. 
Moreover
these time averages are given by a finite number
of physical absolutely continuous
probability measures --- under mild conditions.
This result holds for a fixed noise level, so the question
naturally arises of what can we say when the noise
level shrinks to zero.

We show here (Theorem~\ref{thmA}) that a dynamical
system given by a diffeomorphism $f:M\mapto M$, exhibiting
several (infinitely many) attractors
$\{ \La_k \}_{k=1}^N$, $N\in\natur\cup\{\infty\}$,
has physical probability measures under random perturbations that are
close to invariant probability measures supported on $\La_k$ when
the perturbed maps are close to $f$ if
(i)
\emph{the union of the basins of attraction $\cup_k W^s(\La_k)$ covers
  Lebesgue almost all points of $M$} and
(ii)\emph{ each attractor $\La_k$ is stochastically stable}.
Condition (ii) means that each $\La_k$ supports an $f$-invariant
probability measure $\mu_k$ such that the time averages of continuous
functions $\varphi$, along randomly perturbed orbits \emph{in a
  neighborhood} of the attractor, converge to $\int \varphi \, d\mu_k$
when the perturbed maps are very close to the unperturbed map $f$.
We do not assume that the attractors support $SRB$ measures.

In particular we obtain (Corollary~\ref{corollary})
for the classical Axiom A case
a new result, about the convergence of the distribution of the mean
sojourn time  of randomly perturbed orbits to a specific linear
combination of the $SRB$ measures of the hyperbolic
attractors, with weights  given by the Lebesgue measure of the
basins of attraction of each hyperbolic attractor, when the
perturbations are close to the unperturbed map $f$.

In Theorem~\ref{thmB} we suppose that the number $N$ of attractors is finite
and
replace condition (i) by a weaker condition on
the \emph{trapping regions} of $f$. If for all open sets $U$
satisfying $f(\ov{U})\subset U$ it holds that
$$
\bigcap_{n\ge0} f^n(\ov{U}) \;\;\mbox{contains}\;\;\La_i\;\;
\mbox{for some}\;\; i\in\{1,\ldots, N\},
$$
then we show that the distribution of the mean sojourn times
of randomly perturbed orbits in subsets of $M$ tends to a
linear convex combination of the probability measures
$\mu_1,\ldots,\mu_N$, when the perturbed maps get close to  $f$.

These results favor a global notion of stochastic stability, already
 proposed by e.g. Viana~\cite{V}, as
 opposed to the usual notion of stochastic stability built on single
 attractors:  stochastically stable systems have time averages given by
 measures that are
 close to a  linear convex combination of invariant measures for the
unperturbed system. Moreover this notion of stochastic stability
 extends to settings where infinitely many attractors coexist and even
 to a setting where time averages for the unperturbed map $f$ do not
 exist (cf. Example 3 in subsection~\ref{Bowen}).
Recently the author together with Alves~\cite{AA} proved that some
nonuniformly expanding systems introduced in~\cite{Vn} and~\cite{ABV}
 are stochastically stable in this global sense.

We hope that these results will foster the understanding of the
so called \emph{Newhouse phenomenon} (the coexistence
of infinitely many sinks when generically unfolding
a quadratic homoclinic tangency by a one-parameter family),
 about which little is known apart from its existence.
Taking a one-parameter family $\{f_t\}$ of surface
diffeomorphisms generically unfolding a quadratic
homoclinic tangency at $t=0$, it is known
(cf. Palis-Takens~\cite[Chpt. 6]{PT})
that for generic $t$ close to
$0$ (Baire category) the map $f_t$ exhibits
infinitely many periodic attractors (sinks) whose orbits all pass
through a fixed neighborhood $U$ of a homoclinic tangency point
of $f_0$. In this setting we pose the following

\medskip
\noindent
{\bf Problem:}
\emph{Is there some parameter value $t$ for which
the return map under $f_t$ to a neighborhood $U$ of a quadratic
homoclinic tangency
satisfies the conditions of Theorem~\ref{thmA}?}

\medskip
\noindent
Moreover we lack a characterization of the $f$ invariant measures that
appear as weak$^*$ limits of physical measures when the noise level $\ep$
tends to zero. This is a deep problem in dynamics whose understanding
should enable us to state results similar to Theorems~\ref{thmA}
and~\ref{thmB} in some fairly general setting without assuming the
stochastic stability of the attractors.

\medskip
\noindent
{\bf Problem:}
\emph{Is there some characterization for the $f$ invariant measures obtained
as weak$^*$ limits or weak$^*$ accumulation points of physical
measures when the noise level tends to zero?}

\medskip
\noindent
The necessary definitions and precise statements
of the results
can be found in section~\ref{statement}.
In subsection~\ref{examples} applications of
these results to some simple systems are shown.
Sections~\ref{brpo} and~\ref{stable}
demonstrate Theorems~\ref{thmA} and~\ref{thmB}.

\medskip

\noindent
{\bf Acknowledgments:\/}
This paper benefitted from talks with
Marcelo Viana at IMPA (Rio de Janeiro)
and Michael Benedicks at KTH (Stockholm)
while on leave from Centro de Matem\'atica
da Universidade do Porto (CMUP), Portugal.
These three institutions
also provided partial financial support,
which made this work possible.
Many thanks to Jos\'e Alves, my
colleague at CMUP, for invaluable discussions
and also to the referees whose suggestions helped me to greatly
improve the present text.


\section{Statement of the results}
\label{statement}

\subsection{Physical random perturbations}
\label{physicalpert}

Let $f:M\mapto M$ be a diffeomorphism
of a compact finite dimensional boundaryless
manifold $M$. We suppose $M$ to be endowed
with a smooth Riemannian metric
inducing a distance \emph{dist} and the associated
smooth normalized Riemannian volume $m$,
fixed once and for all.

The random perturbation of $f$ to be considered
is a finite dimensional smooth parameterized family
of diffeomorphisms  given by a $C^1$-map
$F:M\times B_1 \mapto M$,
$(x,t)\mapsto f_t(x)=f(x,t)$,
where $B_1=\{ x\in{\real}^n: \| x \|_2 <1 \}$
is the unit ball on ${\real}^n$
($n$ is not related to $\dimns (M)$)
and for every $t\in B_1$ the map
$f_t:M\mapto M$, $x\mapsto f(x,t)$ is a
diffeomorphism.
This family must contain $f$ at the origin, i.e.
\begin{enumerate}
\item[I)] $f_0\equiv f$.
\end{enumerate}
Moreover considering for a given $\epsilon\in]0,1[$
the \emph{perturbation space}
$$
\De_{\epsilon}=\big\{ \un t = (t_j)_{j\ge1}: t_j\in{\real}^n
\qand \|t_j\|_2\le\ep \big\}
$$
and defining the \emph{perturbed iterate}
of $x\in M$ by the \emph{perturbation vector}
$\subl{t}\in\Delta_{\epsilon}$
as
$$
f^j_{\subl{t}}(x) = f^j(x,\subl{t})=
f_{t_j}\circ \cdots \circ f_{t_1}(x), \qquad j\ge1
\qquad \mbox{and} \quad f^0(x,\subl{t})=x
$$
we impose another condition.
\begin{enumerate}
\item[II)] There are $K=K(\ep)\in \natur$ and
$\xi=\xi(\ep)>0$
such that for all $j\ge K$ it holds
$$
f^j(x,\Delta_{\epsilon}) =
\{ f^j_{\subl{t}}(x) : \subl{t}\in\Delta_{\epsilon} \}
\supset B( f^j(x),\xi ).
$$
\end{enumerate}
This assumption demands that the set of all
perturbed iterates of any given $x\in M$ covers
a full neighborhood of the unperturbed iterate of
$f$, at least after some threshold.
Further, to be able to measure \emph{typical behavior}
we take the infinite product probability
$\nu^{\infty}_{\epsilon}$ on $\Delta_{\epsilon}$,
where $\nu_{\epsilon}={\rm Leb}_n|\ov{B_\ep}$
(the normalized $n$-dimensional Lebesgue volume
measure restricted to the closure of
$B_\ep=\{x\in{\real}^\natur:\|x\|_2<\ep\}$),
and assume the following \emph{non degeneracy} condition
on the map $F$.
\begin{enumerate}
\item[III)] For every $x\in M$ and for all $k\ge K$ we
have $f^k(x,\nu_{\epsilon}^\infty)\ll m$,
\end{enumerate}
where $f^k(x,\nu_{\epsilon}^\infty)$ is the
measure that integrates continuous functions
$\varphi: M \mapto \real$ as
$f^k(x,\nu_{\epsilon}^\infty)\varphi =
\int \varphi( f^k(x,\subl{t}) ) \,
d\nu_{\epsilon}^\infty(\subl{t})$
--- the push-forward of $\nu_{\epsilon}^\infty$
to $M$ by 
$f^k(x,\cdot) : \Delta_{\epsilon} \mapto M$,
a continuous function with respect to the product
topology on $\Delta_{\epsilon}$ and
continuously
differentiable
on each coordinate of $\Delta_{\epsilon}$.
This says that sets of perturbation vectors of
positive $\nu_{\epsilon}^{\infty}$-measure
must send points of $M$ onto positive
volume subsets of $M$ (positive $m$-measure)
after a finite number of iterates.

A family $F$ satisfying items
I, II and III  for all
$\epsilon\in]0,\epsilon_0[$, for some given
fixed $\epsilon_0\in]0,1[$, will be named
a \emph{physical random perturbation of} $f$.
When we fix the value of $\epsilon$ on some
discussion we will just write
\emph{the perturbation of level} $\epsilon$.

As shown in \cite[Examples 1 and 2]{Ar} we can always build
a physical random perturbation of any given
diffeomorphism $f: M\mapto M$ if we
allow a sufficiently big (finite) number
of parameters on the family $F$.

\begin{remark}\label{openmap}
In what follows the arguments will be
written out thinking in terms of
families of diffeomorphisms. However
they hold for continuous families $(f_t)_{t\in B_1}$
of proper continuous maps such that $(f_t)_* m\ll m$
for $\nu$ almost all $t\in B_1$.
\end{remark}

\subsection{Stochastically Stable Attractors}
\label{ssattract}

Let $f: M \mapto M$ be a diffeomorphism
of a compact finite dimensional manifold as
in the previous sections. We will
consider stochastically stable attractors, according
to the following definitions.

\begin{definition}\label{attractor}
An \emph{attractor for} $f$ is a compact  set $\Lambda$ such that
(a) $f(\La)=\La$;
(b) there is a neighborhood $U_{\La}$ of $\La$
such that $f( \tra{U}_{\La} ) \subset U_{\La}$
(a \emph{trapping region});
(c) $\La=\cap_{k\ge1} f^k(\ov{U}_{\La})$ and
(d) there is $x\in\La$
such that ${\rm closure}\{ f^k(x): k\ge 1 \} = \La$.
Moreover the \emph{basin of attraction} of $\La$
is the set
$
W^s(\La) = \{ x\in M: \lim_{k\to+\infty}{\rm dist\,}( f^k(x) , \La) =0\}
$
and we remark that $W^s(\La)$ is always open.
\end{definition}
The following useful property enables us to
\emph{localize} physical probability measures
--- we postpone the proof
to subsection~\ref{localizing}.

\begin{proposition}\label{uniqphysical}
Let $\La$ be an attractor with respect to $f$
and let us take a physical random perturbation
$F$ of $f$.
Then there are an open set $U=U_{\La}$ and $\ep_0\in]0,1[$ such that for
$\ep\in]0,\ep_0[$ we have that
(a) $U$ is
\emph{completely forward invariant}, i.e.
$f_t(\tra{U})\subset U$ for all $t\in \ov{B}_\ep$;
(b) there is only one probability measure $\mu^{\ep}$
such that $\La\subset\supp\mu^{\ep}\subset U$
and for all $x\in U$ and $\nu_{\ep}^\infty$ almost every
$\subl{t}\in\De_{\ep}$ we have for each continuous
$\varphi:M\mapto\real$ that
$
\lim_{n\to\infty} \frac1n \sum_{j=0}^{n-1}
\varphi( f^j(x,\subl{t}) ) = \int \varphi \, d\mu^{\ep};
$
and also
(c) ${\rm Hd}(\supp\mu^{\ep} , \La ) \mapto 0$
when $\ep\mapto 0^{+}$, where $\rm Hd$ is the Hausdorff distance
between compact subsets of $M$.
\end{proposition}

\noindent
We will need to assume the convergence of the
physical probability measures in the
following sense.

\begin{definition}\label{stability}
An attractor $\La$ with respect to $f$ is
\emph{stochastically stable} with respect to a perturbation $F$
if it supports an $f$-invariant probability measure $\mu$
($\supp\mu=\La$) such that 
$$
\mu^{\epsilon} \mapto \mu \;\;\mbox{in the weak}^{*}
\;\mbox{topology when}\;\; \ep\mapto 0^{+},
$$
where the $(\mu^{\ep})_{\ep>0}$ are given by
Proposition~\ref{uniqphysical}.
A stochastically stable attractor will be written
as a pair $(\La,\mu)$.
\end{definition}

\subsection{Infinitely Many Attractors}
\label{inftyattractors}

Let $f:M\to M$ be a diffeomorphism and $F$ be a physical random
perturbation of $f$.
For maps $f$ with several (infinitely many) attractors
we study the distribution of the mean sojourn times of the
randomly perturbed
orbits of a point $x$ with respect to the probability measure
$\nu_\ep^\infty$ on the perturbation vectors $\un t\in\De_\ep$. We note
that for a measurable subset $A\subset M$ the measure
$
\nu_\ep^\infty\left\{ \un t\in\De_\ep : f^j(x,\un t)\in A \right\}
$
is the probability of choosing a perturbation vector that sends $x$
into the set $A$ after $j$ iterates. Then the average
$$
\mu^\ep_n(x)A=\frac1n\sum_{j=0}^{n-1}
\nu_\ep^\infty\left\{ \un t\in\De_\ep : f^j(x,\un t)\in A \right\}
$$
gives
the $\nu_\ep^\infty$-averaged frequency of visits to the set $A$
of the perturbed
orbits of $x$ in $n$ iterates, which defines a Borel probability
measure $\mu_n^\ep$ on $M$.
We define the \emph{mean sojourn time}
of the random orbits of a point $x$ as the probability measure
\begin{equation}
  \label{meansojournforx}
  \mu^\ep(x) =w^*-\lim_{n\to+\infty} \mu^\ep_n(x)
\end{equation}
if the weak$^*$ limit exists.
Likewise we define for a measurable subset $A\subset M$ the
$\nu_\ep^\infty$-averaged
frequency of visits to $A$ of the perturbed orbits in $n$ iterates by
$$
\mu^\ep_n A=\frac1n\sum_{j=0}^{n-1}( m\times\nu_\ep^\infty)
\left\{ (x,\un t)\in M\times\De_\ep : f^j(x,\un t)\in A \right\}
$$
which is a probability measure on $M$ and then define the \emph{mean
  sojourn time} of the random orbits of the system as the probability
measure
\begin{equation}
  \label{meansojourn}
  \mu^\ep =w^*-\lim_{n\to+\infty} \mu^\ep_n
\end{equation}
if the weak$^*$ limit exists.
Since $\mu_n^\ep=\int \mu_n^\ep(x) \,dm(x)$ if
$\mu^\ep(x)=\lim_{n\to\infty}\mu_n^\ep(x)$ exists $m$-a.e., then
$\mu^\ep=\lim_{n\to\infty}\mu_n^\ep=\int \mu^\ep(x) \, dmu(x)$.
When the noise level $\ep$ tends
to zero we have the following properties.

\medskip

\begin{maintheorem}\label{thmA}
Let $f:M\mapto M$ be a diffeomorphism
of a compact finite dimensional boundaryless manifold $M$.
We suppose that a physical random perturbation $F$ of $f$
is given and that there is a  family
of pairwise disjoint compacts $\{ \La_i \}_{i=1}^N$,
$N\in \natur\cup\{\infty\}$, satisfying
\begin{nlista}{\Letra}
\item each $\La_i$ is an attractor
($i=1,\ldots,N$) and
the \emph{basins of attraction} of the
$\La_i$ cover $m$ almost every point of $M$:
$
\bigcup_{i=1}^N W^s(\La_i) = M, \;\; m\;\mbox{mod}\;0.
$
\item each attractor $\La_i$ supports
an $f$-invariant probability measure $\mu_i$ which is
stochastically stable with respect to $F$,
$i=1,\ldots, N$, $N\in\natur\cup\{\infty\}$,
according to definition~\ref{stability}.
\end{nlista}
Then it holds that
\begin{enumerate}
\item for $m$ almost all $x\in M$ there are $\ep_x>0$
and $i=i(x)\in\{1,\ldots,N\}$ such that
for every $\ep\in]0,\ep_x[$ 
the mean sojourn time $\mu^\ep(x)$
of the perturbed orbits of $x$ exists and satisfies
\begin{equation}\label{sojournforxconv}
\mu^\ep(x)=\mu_i^\ep
\;\;\mbox{and}\;\;
\mu^{\ep}(x)\mapto \mu_i
\;\;\mbox{when}\;\; \ep\mapto 0^{+},
\end{equation}
where the convergence is in the
weak$^*$ topology and $\mu_i$ is defined by $\La_i$ as in
Proposition~\ref{uniqphysical}. 
\item for every $\ep>0$ there are
finitely many physical measures
$\mu_1^\ep,\dots,\mu_l^\ep$,
$l=l(\ep)$
for which
there exist $\be_1^\ep,\dots,\be_l^\ep\ge0$
with
$\be_1^\ep +\cdots + \be_l^\ep = 1$ such that the mean sojourn time
for the randomly perturbed system $\mu^\ep$ exists and satisfies
\begin{equation}\label{meansojourn=linear}
\mu^\ep=
\be_1^\ep\mu_1^\ep+\cdots+
\be_l^\ep\mu_l^\ep
\end{equation}
and
\begin{equation}\label{meansojournconv}
\lim_{\ep\to0} \mu^\ep =
\sum_{i=1}^N
m( W^s(\La_i) )\cdot \mu_i,
\end{equation}
where the convergence is always
in the weak$^*$ topology.
\end{enumerate}
\end{maintheorem}

\noindent
In this setting $\mu^\ep(x)$  is  near a
$f$-invariant probability measure $\mu_i$ for
$m$ almost every $x\in M$ when the
noise level is close to zero.
Moreover if $\vfi:M\to\real$ is continuous, then
$\mu^\ep(x)\vfi$ is the mean time average of $\varphi$
along almost all random orbits
of $x$ and thus we know that for $m$ almost all $x\in M$
these averages are near
$\int \varphi\, d\mu_i$
for some $i=i(x)\in\{1,\ldots,N\}$
when the noise level is close to zero.
In addition \refer{meansojourn=linear} is saying
that $\mu^\ep$ is  a  linear convex combination of the physical measures of
the random system, and
according to
\refer{meansojournconv}, when the noise level
is close to zero, $\mu^\ep$ tends
to a linear convex combination
of $f$-invariant measures, supported
on the attractors, whose weights are
given by the volume of the respective
stable sets.

Since an Axiom A system satisfies both Conditions $\un A$ and $\un B$
with a finite number of hyperbolic
attractors, each of which supporting a unique $SRB$ measure which is
stochastically stable,
we have the following corollary.

\begin{corollary}\label{corollary}
Let $f:M\to M$ be a diffeomorphism satisfying Axiom A and let
$(\La_1,\mu_1)$, ..., $(\La_N,\mu_N)$ be the
hyperbolic attractors for
$f$ and their
respective $SRB$ measures. Then we have that the mean sojourn time for
the randomly perturbed system by any physical random perturbation $F$
satisfies
$$
\lim_{\ep\to0} \mu^\ep = m(W^s(\La_1))\cdot\mu_1+\ldots+
m(W^s(\La_N))\cdot\mu_N
$$
where the limit is in the weak$^*$ topology.
\end{corollary}

\noindent
Now we replace Condition $\un A$ by a weaker condition and suppose
that the number of attractors is finite.

\begin{maintheorem}\label{thmB}
Let $f:M\mapto M$ be a diffeomorphism
with a finite family of
stochastically stable attractors
$(\La_1,\mu_1),\ldots,(\La_N,\mu_N)$,  with respect to
a given physical random perturbation
$F$ (and $f$-invariant probability measures
$\mu_1,\ldots,\mu_N$), such that
\begin{itemize}
\item[$\subl C$.] given any open subset $U\subset M$, $U\neq M$,
satisfying $f(\tra{U})\subset U$, it holds that
$$
\bigcap_{k\ge0} f^k (U)\;\;\mbox{contains}\;\; \La_i\;\;
\mbox{for some}\;\; i\in\{1,\ldots,N\}.
$$
\end{itemize}
Then for every $x\in M$ the mean sojourn time
$\mu^{\ep}(x)$ satisfies
\begin{equation}\label{thmB1}
\mu^{\ep}(x) \mapto \mbox{convex hull}\;
\{\mu_{i_1},\ldots,\mu_{i_k}\}
\;\;\mbox{in the weak}^{*}\;\mbox{sense, when}\;\; \ep\mapto 0^{+}.
\end{equation}
\end{maintheorem}

\medskip
\noindent
The meaning of the convergence to the
convex hull of the measures is to say that the
integrals $\int \vfi\,d\mu^\ep(x)$
are close to the set of
linear convex combinations
$$
\left\{\al_1\int\vfi\,\mu_{i_1}+\cdots+
\al_k\int\vfi\,\mu_{i_k},\;
\al_1,\ldots,\al_k\ge0,\;
\al_1+\cdots+\al_k=1\right\}
$$
for all continuous $\vfi:M\mapto\real$
and small enough $\ep>0$.


\subsection{Examples}
\label{examples}

\begin{figure}[h]
 
\centerline{\psfig{figure=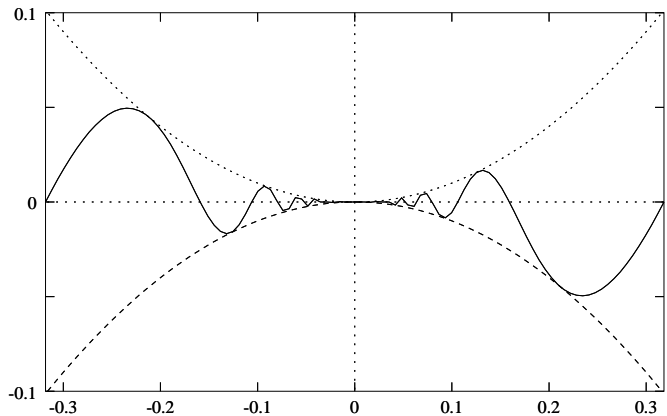,width=3.5in}}

\caption{\label{sinus} A map whose gradient flow has infinitely many sinks} 
\end{figure}
This example illustrates with a simple map the meaning of
Theorems~\ref{thmA} and~\ref{thmB}.

\begin{example}\label{circlemap}
We take the function
$\phi:[-\frac1{\pi},\frac1{\pi}]\mapto \real$,
$s\mapsto s^4\sin(\frac1{s})$
and identify the endpoints of the interval to get a $C^1$ map
$\vfi:S^1\mapto \real$ (cf. Figure~\ref{sinus}).
Letting $\nabla\vfi:S^1\mapto\real$ be the gradient of $\vfi$
and considering the \emph{gradient flow} given by
$\dot{x}=\nabla\vfi(x)$, we easily see the 
\emph{time one map} $f:S^1\mapto S^1$ of this
flow to have infinitely many sources and sinks.
In fact, these correspond to the local maxima and
local minima
of $\phi$, respectively 
(the zeros of $\nabla\vfi$ except the
point $p$ in $S^1$ corresponding to
$0$ in $[-\frac1{\pi},\frac1{\pi}]$).
Moreover, the
basin of attraction of the local minima of $\phi$
extends precisely until the next neighboring local
maxima, and thus the union of all these basins covers
$S^1$ with the exception of the denumerable set of
local maxima and the point $p$.

We conclude that we are in the conditions of Theorem~\ref{thmA}
since every sink is a stochastically
stable attractor.
Hence this system $f$ under the physical random
perturbation given by the $C^1$ family
$F:S^1\times[-1,1]\mapto S^1$, $(x,t)\mapto R_t(f(x))$
is stochastically stable globally as in Theorem~\ref{thmA},
where $R_t:S^1\mapto S^1$ is the rotation of
angle $t$.

We remark that the basins of attraction of the sinks nearer to the
origin do not persist under random perturbations. However for a given
sink, when the level of noise is below a certain threshold, an
invariant region under random perturbations must appear containing a
neighborhood of the sink and supporting a physical probability
measure. As the level of noise is reduced these invariant regions
occupy a larger and larger proportion of the space and their physical
measures tend to the Dirac delta probability measures concentrated on
the sinks. Moreover the weak$^*$ limit of the mean sojourn times of
all points $x$ when $\ep\to0$ cannot give positive weight to more than
two neighboring sinks, with the exception of $x=0$ for which
$\mu^\ep(0)$ must be supported in a neighborhood of $0$ whose diameter
shrinks to zero when $\ep\to0$. Hence $\mu^\ep(0)\to\de_0$ when
$\ep\to0$ in the weak$^*$ topology.
\end{example}

\begin{figure}[h]
 
\centerline{\psfig{figure=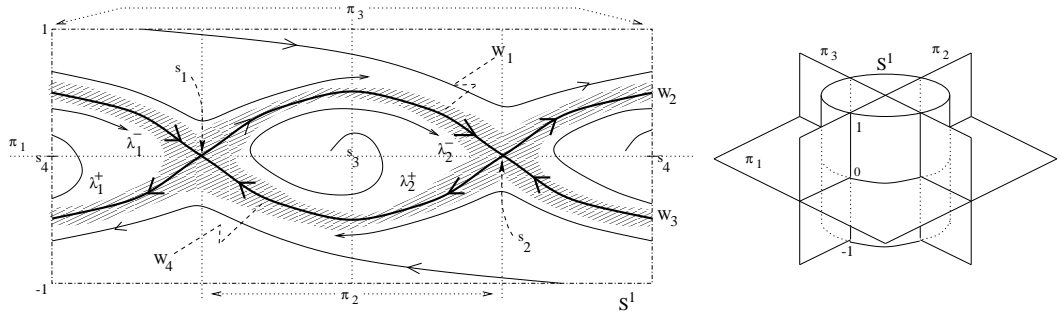,width=5.5in}}

\caption{\label{Bowenfig} A sketch of the flow example of Bowen
with symmetries}
\end{figure}
\medskip

\noindent
The following example shows how the global notion of stochastically
stable systems is flexible enough to encompass systems for which time
averages are \emph{not defined} for Lebesgue almost every point.

\begin{example}\label{Bowen}
The flow $(\phi_t)_{t\ge0}$
sketched on the left hand side of
figure~\ref{Bowenfig} is
attributed to Bowen
(cf.~\cite{Tk}) and is
well known for showing that
\emph{Birkhoff averages may not exist almost
everywhere}.
In this system time averages exist only for
the sources $s_3, s_4$ and for the set of separatrixes
and saddle equilibria
$W=W_1\cup W_2 \cup W_3 \cup W_4 \cup \{s_1,s_2\}.$
Moreover the orbit $(\phi_t(x))_{t\ge0}$ of
each $x$ not in $W$ and different from $s_3, s_4$
tends to $W$ as $t\to +\infty$.

Letting $f=\phi_1$ denote the time $1$ map of the
flow, we see that $W$ is a compact $f$-invariant
set and $W=\cap_{n\ge0} f^n(U)$
 for all sufficiently small
neighborhoods $U$ of $W$.
Hence $W$ lacks only transitiveness to become
an attractor according to definition~\ref{attractor}.
However, defining
$F:(S^1\times [-1,1])\times{\real}^2\longrightarrow
S^1\times [-1,1]$ by
$$
(x,y,(v_1,v_2))\mapsto (x+v_1 \bmod 1, y+v_2)
$$
where we identify $S^1$ with $\real/\relativ$,
and letting
$B_\ep=\{ (v_1,v_2)\in{\real}^2:v_1^2+v_2^2\le\ep^2\}$
for small enough $\ep>0$, we have that
the family $F$ is a physical
random perturbation of $f$.
In addition, as shown in~\cite[Prop. 10.1]{Ar},
there is a \emph{unique} physical probability
measure $\mu^\ep$ for the perturbed system whose
support $\supp(\mu^\ep)$ is a neighborhood of $W$,
sketched in Figure~\ref{Bowenfig}
as a shaded region.
This neighborhood is small for small $\ep>0$,
satisfying all items of Proposition~\ref{uniqphysical},
although $W$ is not a transitive set.

The reason for this is that $W$ is
\emph{chain-transitive} or
\emph{transitive with respect to pseudo-orbits}.
Indeed for all $\de>0$ and every pair
$x,y\in W$ there is a $\de$-pseudo-orbit
$(z_k)_{k=0}^n$, for some $n\in\natur$,
such that $z_0=x$ and $z_n=y$,
i.e., every pair of points of $W$ can
be joined by some $\de$-pseudo-orbit.
Since random orbits are $\de$-pseudo-orbits
for some $\de>0$, it is easily checked
that the arguments of Proposition~\ref{uniqphysical}
can be repeated in this setting.

Finally we observe, on the one hand, that
every weak$^*$ accumulation point $\mu$ of
$(\mu^\ep)_{\ep>0}$, when $\ep\to 0$, is
a $f$-invariant probability measure and,
on the other hand, that $\mu$ must be a
convex linear combination
$\la\de_{s_1}+(1-\la)\de_{s_2}$, $0\le\la\le1$,
of the Dirac measures on the saddle fixed points.
In fact every $f$-invariant probability measure
must be a convex linear combination of
$\de_{s_1},\de_{s_2},\de_{s_3},\de_{s_4}$
and $\supp(\mu)\subset\supp(\mu^\ep)$
for all small $\ep>0$.

This means that $\de_{s_1},\de_{s_2}$ are
\emph{stochastically stable}! Indeed
$C=\{ \la\de_{s_1}+(1-\la)\de_{s_2} : 0\le\la\le1\}$
is a closed subset of $\mathcal{P}(M)$,
the compact metrizable space
of all probability measures on
$M=S^1\times[-1,1]$ with the weak$^*$ topology,
and 
$(\mu^\ep)_{\ep>0}$ is a family whose
accumulation points all lie in $C$.
Hence by compactness of $C$ and
because the weak$^{*}$ topology
on $\mathcal{P}(M)$ is a metrizable
one
$$
\mu^\ep \longrightarrow
\mbox{ convex hull }\{ \de_{s_1},\de_{s_2} \}
\qquad \mbox{when   } \ep\to 0,
$$
meaning that for all continuous
$\vfi:M\to\real$
and  for small enough $\ep>0$
$$
\int \vfi\,d\mu^\ep \qquad
\mbox{is close to} \qquad
\Big\{
\la\vfi(s_1)+(1-\la)\vfi(s_2) : 0\le \la\le1
\Big\}.
$$
\end{example}


\section{Behavior of randomly perturbed orbits}
\label{brpo}

We explain here the main ingredients involved in the proof of
Theorems~\ref{thmA} and~\ref{thmB}
recalling the main lemmas used in the proof of a general result
in~\cite{Ar}.
This result is stated below and guarantees the existence of a
finite number of physical measures for general physical random perturbations.
In  subsection~\ref{localizing}
we present a proof of proposition~\ref{uniqphysical}

\begin{theorem}\label{thm0}
Under a physical random perturbation
of fixed level $\epsilon>0$ of some diffeomorphism
$f: M\mapto M$, there are a finite number
$\mu_1,\ldots,\mu_l$ of \emph{physical} probability measures in
$M$ such that
\begin{enumerate}
\item $\supp \mu_i \cap \supp \mu_j = \emptyset$,
$1\le i < j \le l$;
\item $\mu_i\ll m$ and
$\supp (\mu_i)$ contains a ball of radius $\xi=\xi(\ep)>0$, $i=1,\ldots,l$;
\item for every $x\in M$ there are open sets
$V_1=V_1(x),\ldots,V_l=V_l(x)\subset \Delta_{\epsilon}$
such that
\begin{enumerate}
\item $V_i\cap V_j = \emptyset$, $1\le i < j \le l$;
\item $\nu_{\epsilon}^\infty (\Delta_\epsilon \setminus
(V_1\cup\ldots\cup V_l))=0$;
\item for every $1\le i\le l$ and all $\subl{t}\in V_i$
we have for every continuous $\varphi:M\mapto\real$
\begin{equation}\label{orbaverage}
\frac1n\sum_{j=0}^{n-1} \varphi( f^j(x,\subl{t})) \mapto
\int \varphi \, d\mu_i \;\;\mbox{when}\;\;
n\to\infty;
\end{equation}
\item $x\in\supp\mu_i \impl V_i(x)=\De_{\ep}$,
$\nu_{\ep}^{\infty}$-mod $0$, for $i=1,\ldots,l$;
\item the $V_i(x)$ depend continuously on $x$
with respect to the metric
$d_{\ep}(A,B)=\nu_{\epsilon}^{\infty}(A \difsim B)$,
$A$, $B$ measurable subsets of $\Delta_{\epsilon}$.
\end{enumerate}
\end{enumerate}
Moreover  every completely forward invariant nonempty set $C$
($f_t(C)\subset C$
for all $t\in B_{\ep}$) contains $\supp\mu_i$ for
some $i\in\{1,\ldots,l\}$ and
$\SM_i=\supp\mu_i$ is a
\emph{trapping region} with respect to
$f^K\equiv f_0^K$, i.e.
$f_t^K( \SM_i ) \subset \inter{ \SM_i }$
for all $i=1,\ldots,l$
(where $K=K(\ep)$ is given by condition
II on $F$).
\end{theorem}

\begin{proof}
See~\cite{Ar}.
\end{proof}

\medskip
\noindent
In what follows we state the main lemmas used to prove
Theorem~\ref{thm0} in~\cite{Ar}.

\subsection{Decomposing invariant probability measures}
\label{building}

We assume that a diffeomorphism $f$
and
a random perturbation $F$ of $f$ with
fixed noise level $\ep>0$ are given.
We set $S_\ep:M\times\De_\ep\to M\times\De_\ep$ to be the skew-product
map
$(x,\un{t})\mapsto (f_{t_1}(x),\sigma \un{t})$
where $\sigma$ is the usual left shift on sequences $\un
t=(t_1,t_2,\ldots)$.

\begin{proposition}
\label{asymptotic}
Let $D$ be a completely forward invariant subset of $M$
($D=M$ is allowed).
Given any $x\in D$ every $\mbox{weak}^{*}$ accumulation
point $\mu$ of the sequence
\begin{equation}\label{averageofaverage}
\left( k^{-1} \sum_{j=0}^{k-1} f^j(x,\nu_{\ep}^{\infty})
\right)_{k=1}^{\infty}
\end{equation}
is a probability measure such that
$\mu\times\nu_{\ep}^{\infty}$ is $S_{\ep}$-invariant
and $\SM=\supp\mu$ is completely forward invariant, a trapping region for
$f^K\equiv f^K_0$ and $\SM\subset D$.
\end{proposition}

\begin{proof}
V.~\cite[sects. 6, 7]{Ar}.
The $S_\ep$-invariance is trivial by construction of $\mu$
and definition of $S_{\ep}$. The complete forward invariance of $\SM$ is an
easy consequence of the form of $S_\ep$.
If we take $K$ given by the property II of a physical
random perturbation, that is
$f^K(x,\De_{\ep})\supset B(f^K(x),\xi)$
for all $x\in D$, then we easily see that $\SM$ is a trapping region
for $f^K$.
\end{proof}


\begin{definition}
\label{ergodicbasin}
Given a probability measure $\mu$ such
that $\mu\times\nu_{\ep}^{\infty}$
is $S_{\ep}$-invariant, we say that the set
$
B(\mu)=\left\{
x\in M : k^{-1} \sum_{j=0}^{k-1}
\de_{f^j(x,\subl{t})} \stackrel{w^{*}}{\mapto} \mu
\;\;\mbox{when}\;\; k\to\infty
\;\;\mbox{for}\;\; \nu_{\ep}^{\infty}\mbox{a.e.}\;
\subl{t}\in B_{\ep}
\right\}
$ is the basin of $\mu$.
\end{definition}

\begin{proposition}
\label{decomposition}
There are a finite number $\mu_1,\ldots,\mu_l$
of probability measures in $M$ such that
for all $1\le i,j \le l$
\begin{enumerate}
\item $\mu_i\times\nu_{\ep}^{\infty}$
is $S_{\ep}$-ergodic, $\supp\mu_i \cap \supp\mu_j = \emptyset$ and
$\mu_i(\pa \supp\mu_i)=0$;
\item $B(\mu_i)\supset \supp\mu_i$;
\item for every  probability
measure $\mu$ such that $\mu\times\nu_\ep^\infty$ is
$S_{\ep}$-invariant there are
$\al_1,\ldots,\al_l\ge0$ with
$\al_1+\ldots+\al_l=1$ such that
$\mu=\al_1\mu_1 + \ldots + \al_l\mu_l$.
\end{enumerate}
\end{proposition}

\begin{proof}
See~\cite[Sect.8]{Ar}.
\end{proof}

\begin{remark}
\label{V=Delta}
We observe that item 2 above corresponds to
item 3d of Theorem~\ref{thm0}.
\end{remark}
The coefficients of item 3 of last
proposition have a natural meaning.

\begin{lemma}\label{coefficients}
Let $\mu(x)$ be a $\mbox{weak}^{*}$
accumulation point of the sequence~\refer{averageofaverage}.
Then
$$
\al_i = \al^\ep_i(x) =
\nu_{\ep}^{\infty}
\{
\subl{t}\in \De_{\ep} :
\exists k\ge1 \;\;\mbox{such that}\;\;
f^k(x,\subl{t})\in\supp\mu_i
\},
\; i=1,\ldots,l.
$$
\end{lemma}
Consequently because $\al^\ep_i(x)$
depends only on $x$, we see that there is
a single $\mbox{weak}^{*}$ accumulation
point for the sequence~\refer{averageofaverage}
and thus we conclude that
 $\mu^\ep(x)$ is the weak$^*$ limit of the
sequence~\refer{averageofaverage}
when $k\to\infty$.

\begin{proof}
We have
$\mu(x)={\rm w^{*}}-\lim_{j\to\infty}
k_j^{-1} \sum_{j=0}^{k_j-1} f^j(x,\nu_{\ep}^{\infty})$
for some sequence $k_1<k_2<\ldots$ of integers.
According to proposition~\ref{decomposition}
there are $\al_1,\ldots,\al_l\ge0$ with
$\al_1+\ldots+\al_l=1$ and $\mu_1,\ldots,\mu_l$
satisfying $\mu(x)=\al_1\mu_1+\ldots+\al_l\mu_l$.
Since the $\mu_1,\ldots,\mu_l$ have pairwise disjoint
compact supports, it must be that
$\mu_i( \supp\mu_i ) = \de_{ij}$ (Kronecker delta),
$i,j=1,\ldots,l$.
Therefore $\mu(x)(\supp\mu_i)=\al_i$,
$1\le i\le l$.

The construction of $\mu(x)$ tells us then
that, because $\mu_i(\pa\supp\mu_i)=0$,
$$
\al_i=\al_i(x)=\lim_{h\to\infty}
\frac1{k_h} \sum_{j=0}^{k_h-1}
\nu_{\ep}^{\infty}(V_i^j(x))
\;\;\mbox{where}\;\;
V_i^j(x)=
\{
\subl{t}\in\De_{\ep}: f^j(x,\subl{t})\in \supp\mu_i
\}.
$$
Observing that $V_i^j(x)\subset V_i^{j+1}(x)$ by
the invariance of $\supp\mu_i$,
$i=1,\ldots,l; j\ge1$, we see that defining
$
V_i(x)=\{ \subl{t}\in \De_{\ep} : \exists k\ge1
\;\mbox{s.t.}\;
f^k(x,\subl{t})\in\supp\mu_i \} = \cup_{j\ge1} V_i^j(x)
$
we arrive at
$\nu_{\ep}^{\infty}(V_i) =
\lim_{j\to\infty} \nu_{\ep}^{\infty} (V^j_i(x))
=\lim_{k\to\infty} k^{-1} \sum_{j=0}^{k-1}
\nu_{\ep}^{\infty}( V_i^j(x)) = \al_i(x)$.
\end{proof}

\medskip

Finally, after these propositions,
for any given $x\in M$ the
\emph{unique probability} $\mu(x)$
decomposes in a \emph{unique way}
$\al_1\mu_1+\ldots+\al_l\mu_l$ and
$\De_{\ep}$ has  a $\nu_{\ep}^{\infty}$ mod $0$
partition $V_1,\ldots,V_l$
since,
clearly, the $V_1=V_1(x),\ldots,V_l=V_l(x)$ are
pairwise disjoint by definition and their total
$\nu_\ep^\infty$ measure is $\al_1+\cdots+\al_l=1$.

\subsection{Localizing physical probability measures}
\label{localizing}
Now we are able to prove Proposition~\ref{uniqphysical}.

\begin{rproof}{\ref{uniqphysical}}
We suppose that $\La$ is an attractor, that for a given
open set $U$ we have
$\La=\cap_{k\ge0} f^k(U)$ and $f(\tra{U})\subset U$
and that $F$ is a physical random perturbation
of $f$.
The parameterized family 
$F$ given by the map
$F:M\times B_1(0)\mapto M$
is such that, by continuity
\begin{equation}\label{moves2zero}
\diame{ f(x,\De_{\ep})} =
\sup \big\{
{\rm dist\,} \left( f_t(x),f_s(x) \right),
t,s\in B_{\ep} \big\} \mapto 0
\;\;\mbox{when}\;\; \ep\mapto 0^{+}.
\end{equation}
Therefore for sufficiently small $\de_0>0$
the open set $U$ is a trapping region
with respect to $f_t$
for all $\ep\in]0,\de_0[$ and every $t\in B_{\ep}$.
In particular this means that $U$ is a completely forward invariant set and
then there is some physical probability
$\mu^{\ep}$ with $\supp\mu^{\ep}\subset U$ by Theorem~\ref{thm0}.

It is easy to see that $\SM^{\ep}=\supp\mu^{\ep}$
is a compact forward $f_t$-invariant set
$f_t(\SM^{\ep})\subset \SM^{\ep}$
for each $t\in B_{\ep}$ and $\ep\in]0,\de_0[$.
Indeed $(\supp\mu^{\ep})\times\De_{\ep}=
\supp(\mu^{\ep}\times\nu_{\ep}^\infty)$ is
$S_{\ep}$-invariant and the definition
of $S_{\ep}$ gives the  invariance of $\SM^\ep$.
In particular $\SM^{\ep}$ is
forward $f$-invariant ($f\equiv f_0$).
On the one hand  $\La\supset \cap_{k\ge0} f^k(\SM^{\ep})$
 since $\La$ is the maximal
forward $f$-invariant set in $U\supset\SM^{\ep}$.
On the other hand
$\SM^{\ep}$ is a $f_t^K$-trapping region
(by Theorem~\ref{thm0})
and so $\cap_{k\ge0} f^k(\SM^{\ep}) \subset \inter{\SM^{\ep}}$.
Thus $\La\cap \inter{\SM^{\ep}} \neq\emptyset$.
Because $\La$ is $f$-transitive and $\SM^{\ep}$
is forward $f$-invariant, it must
be that $\La\subset\inter{\SM^{\ep}}$.

Observing that this holds for every physical
probability measure $\mu^{\ep}$ with $\supp\mu^{\ep}\subset U$,
and that for every random perturbation of fixed
level $\ep$ the supports of the physical probability measures
are pairwise disjoint by item 1 of Theorem~\ref{thm0},
we deduce that there is a
\emph{unique physical probability measure} $\mu^{\ep}$ inside $U$.
This provides items $(a)$ and $(b)$ and the
uniqueness part of Proposition~\ref{uniqphysical}.
In addition the orbital averages of any $x\in U$
must satisfy item $(c)$, since the random orbits
cannot leave $U$ and their time averages exist
and are given by a physical probability measure
(Theorem~\ref{thm0} again).
It must be $\mu^{\ep}$ by the previous uniqueness
arguments.

To get $(d)$ we first note that $\La\subset\supp(\mu^\ep)$ for all
small $\ep>0$. Now we just need to show that for any given
neighborhood $W$ of $\La$ we may find $\de_0$ such that
$\supp(\mu^\ep)\subset W$ for every $\ep\in]0,\de_0[$.
It is easy to see that the set
$U_l=\cap_{k=0}^l f^k(U)$ satisfies
$
\dista{U\setminus U_l}{\La} \mapto 0
\;\;\mbox{when}\;\; l\mapto\infty
$, i.e. $W\supset U_l$ for big enough $l$,
and $\La\subset f(\tra{U}_l)\subset U_l$
for all $l\ge1$.
Hence~\refer{moves2zero} ensures
that for given $l\ge1$ there is $\de_0>0$
such that $f_t(\tra{U}_l)\subset U_l$
for all $t\in B_{\ep}$ and  $\ep\in]0,\de_0[$.
Thus  $\supp\mu_{\ep}\subset U_l$ for all
$\ep\in]0,\de_0[$ with $\de_0(l)\mapto 0$ as
$l\mapto\infty$. This shows that
$\supp(\mu^\ep)\to\La$ in the Hausdorff
topology when $\ep\to0$.
\end{rproof}



\section{Stochastic stability}
\label{stable}

Here we assume that $f$ is a diffeomorphism
satisfying  conditions $\un A$ and $\un B$
of Theorem~\ref{thmA}.

\begin{lemma}\label{basin=stable}
If $\mu_i^\ep$ is the physical probability measure provided by
Proposition~\ref{uniqphysical} for all small enough $\ep>0$, then
$m[ W^s(\La_i) \setminus B(\mu_i^{\ep})] \mapto 0$
when $\ep\to 0^{+}$ for $i\in\{1,\ldots,N\}$.
\end{lemma}

\begin{proof}
We fix $i\in\{1,\ldots,N\}$ and observe that because $\La=\La_i$ is an
attractor there is a trapping region $U$ such that $\La=\cap_{n\ge1}
f^n(\ov{U})$. By definition of basin of attraction
$W^s(\La)=\cup_{j=0}^\infty f^{-j}(\ov{U})$ and
since $U$ is a trapping region
$f^{-j-1}(\ov{U})\supset f^{-j}(\ov{U})$ for all $j\ge0$.
By Proposition~\ref{uniqphysical} there is $\de_0>0$ such that for
every $\ep\in]0,\de_0[$ it exists a physical measure
$\mu^\ep=\mu^\ep_i$ satisfying $\supp(\mu^\ep)\subset U$ and
$B(\mu^\ep)\supset U$. We have to show that for all $j\ge1$ there is
$\de_1\in]0,\de_0[$ such that for every $\ep\in]0,\de_1[$ it holds
that $B(\mu^\ep)\supset f^{-j}(\ov{U})$.

Let us fix an arbitrary integer $j\ge1$. Then for all $x\in
f^{-j}(\ov{U})$ we know that, because $U$ is a trapping region and by
the continuity of the map
$f^{j+1}(x,\cdot):\De_{\de_0}\to M$,
there is $\de_1\in]0,\de_0[$ such that
$f^{j+1}(x,\De_\ep)\subset U$ for all $\ep\in]0,\de_1[$.
By compactness of $f^{-j}(\ov{U})$ ($f$ is a diffeomorphism) the value
$\de_1$ may be taken to depend on $j$ alone.
It is easy to see that if $f^{j+1}(x,\un t)\in B(\mu^\ep)$ for fixed $j\ge1$
and for $\nu_\ep^\infty$ almost all $\un
t\in\De_\ep$, then $x\in B(\mu^\ep)$ also.
We conclude that $B(\mu^\ep)\supset f^{-j-1}(\ov{U})\supset
f^{-j}(\ov{U})$
for all $\ep\in]0,\de_1[$.
\end{proof}

\begin{remark}\label{nondecreasingbasin}
The proof of this lemma  shows that $B(\mu_i^\ep)$ is a
nondecreasing function of $\ep$.
\end{remark}

\noindent
Now for $i\in\{1,\ldots, N\}$ we define the thresholds
\begin{equation}
  \label{thresholds}
  \ep_i= \sup\left\{
\de\in]0,1[: \forall \ep\in]0,\de[ \;\; \exists! \mu_i^\ep \;\;\mbox{such
  that}\;
\La_i\subset\supp(\mu_i^\ep)\subset\supp(\mu_i^\de)
\right\}
\end{equation}
(the set is nonempty by Proposition~\ref{uniqphysical}).
We observe that when the number $N$ of attractors is infinite, it must
be that $\inf_i\{\ep_i\}=0$ for otherwise we would have an infinite
number of  physical measures $(\mu_i^\ep)_{i\ge1}$
for $\ep\in]0,\inf_i\{\ep_i\}[$ such that $(\supp(\mu_i^\ep))_{i\ge1}$
are pairwise disjoint and each contains
a ball of radius $\xi=\xi(\ep)>0$ (by Theorem~\ref{thm0}). This is
impossible in a compact space.
Now we may reindex the attractors in such a way that
$\ep_1\ge\ep_2\ge\ldots$ holds.
If we define $N(\de)=\sup\{i\ge1 \; : \; \ep_i\ge\de, i\ge1\}$
and $\SV_\de=\cup_{i=1}^{N(\de)} B(\mu_i^\de)$, then
$
m(\SV_\de)\to 1$ when $\de\to0$,
because Condition $\un A$ of Theorem~\ref{thmA} and
Lemma~\ref{basin=stable}
easily imply that
$
\sum_{j=1}^{N(\ep)} m\big( B(\mu_i^\ep) \big) \to 1
$
when $\ep\to 0$.

It is now easy to prove item 1 of Theorem~\ref{thmA} because the
weak$^*$ limit of the sequence~\refer{averageofaverage} is the same as
the limit in the definition~\refer{meansojournforx} of mean sojourn
times.
Indeed, according to Lemma~\ref{basin=stable} and
Remark~\ref{nondecreasingbasin}, for $m$ almost every $x\in M$ there is
$\ep_x>0$ such that if $\ep\in]0,\ep_x[$, then $x\in
B(\mu_i^\ep)$ for some $i=i(x)\in\{1,\ldots,N\}$.
Hence $\mu^\ep(x)=\mu_i^\ep$ and Condition $\un B$ says that
$\mu_i^\ep\to\mu_i$ when $\ep\to0$ in the weak$^*$ topology.

For the proof of item 2 of Theorem~\ref{thmA} we observe that the
probability measure $\mu^\ep$ satisfies for all continuous
$\vfi:M\to\real$
\begin{equation}
  \label{limitdoubleint}
  \mu^\ep(\vfi)=\lim_{n\to+\infty}\frac1n \sum_{j=0}^{n-1}
\int \!\!\! \int \vfi( f^j(x,\un t)) \, d\nu_\ep^\infty(\un t) dm(x)
\end{equation}
if the limit exists for fixed $\ep>0$.
On the one hand, for all $x\in M$ we know that
$\lim_{n\to\infty} n^{-1}\sum_{j=0}^{n-1} \vfi(f^j(x, \un
t))\,d\nu_\ep^\infty(\un t) =
\al_1^\ep(x)\mu_1^\ep(\vfi)+\ldots+\al_l^\ep(x)\mu_l^\ep(\vfi)$ by
Proposition~\ref{asymptotic} and item 3 of
Proposition~\ref{decomposition}. It is clear that
$\sup_x\big| n^{-1}\sum_{j=0}^{n-1} \vfi(f^j(x, \un
t))\,d\nu_\ep^\infty(\un t) \big| = \| \vfi \|_0$ for all $n\ge1$
and~\refer{limitdoubleint} equals
$\int[\lim_{n\to\infty} n^{-1}\sum_{j=0}^{n-1} \vfi(f^j(x, \un
t))\,d\nu_\ep^\infty(\un t)]dm(x)$
by Lebesgue's
Dominated Convergence Theorem. So~\refer{limitdoubleint} also
equals $\be_1^\ep\mu_1^\ep(\vfi)+\ldots+\be_l^\ep\mu_l^\ep(\vfi)$ if
we define $\be_i^\ep=\int \al_i^\ep\, dm$, $i=1,\ldots, l=l(\ep)$ ---
note that $\be_i^\ep$ is well defined by item 3e of Theorem~\ref{thm0}.
This proves~\refer{meansojourn=linear}.

On the other hand~\refer{meansojournconv} is easy when $N$ is finite
since $\mu^\ep=
\be_1^\ep\mu_1^\ep+\cdots+
\be_N^\ep\mu_N^\ep$ for $\ep\in]0,\inf\{\ep_1,\ldots,\ep_N\}[$ and
Condition $\un B$ suffices to get the convergence if we show that
$\lim_{\ep\to0} \be_i^\ep=m(W^s(\La_i))$, $i=1,\ldots, N$.
But in fact
$\al_i^\ep(x)=\nu_\ep^\infty(V_i(x))$ by the proof
of Lemma~\ref{coefficients} and $\al_i^\ep(x)=1$ if $x\in
B(\mu_i^\ep)$ by item 3 of Theorem~\ref{thm0}. So $\be_i^\ep\ge
m(B(\mu_i^\ep))$ and then $\liminf_{\ep\to0}\be_i^\ep\ge m(W^s(\La_i))$
by Lemma~\ref{basin=stable}.
Since $m$ is a probability measure
Condition $\un A$ together with Lemma~\ref{basin=stable} imply that
$\limsup_{\ep\to0}\be_i^\ep \le m(W^s(\La_i))$.

Now we assume that $N$ is infinite. It is obvious that $l(\ep)\ge
N(\ep)\to+\infty$ when $\ep\to0$. Let $\be_i=m(W^s(\La_i))$ for
$i\ge1$.  Fixing a continuous
map $\vfi:M\to\real$ we want to bound
$D(\vfi,\ep)=\big|\sum_{i=1}^{l(\ep)} \be_i^\ep\mu_i^\ep(\vfi)-
\sum_{i\ge1} \be_i\mu_i(\vfi) \big|$
by a quantity tending to zero when $\ep\to0$.

Let us fix $\de>0$. Then it exists $I\ge1$
such that $\sum_{i\ge I}\be_i<\de$ by Condition $\un A$ and
\begin{eqnarray*}
{\textstyle
  D(\vfi,\ep)}
& \le & 
  {\textstyle
  \big|
    \sum_{i=1}^I \be_i^\ep\mu_i^\ep(\vfi)-
    \sum_{i=1}^I \be_i\mu_i(\vfi)
  \big| +
  \big|
    \sum_{i=I+1}^{l(\ep)} \be_i^\ep\mu_i^\ep(\vfi)-
    \sum_{i=I+1}^{l(\ep)} \be_i\mu^\ep_i(\vfi)
  \big|
  } \\
& & 
  {\textstyle
    +
    \big|
    \sum_{i=I+1}^{l(\ep)} \be_i\mu_i^\ep(\vfi)-
    \sum_{i>I} \be_i\mu_i(\vfi)
    \big| = D_1 + D_2 + D_3.
  }
\end{eqnarray*}
for $\ep\in]0,\ep(I)[$ where $\ep(I)=\sup\{ \ep\in]0,1[: l(\ep) >
I\}$. By previous arguments $D_1\to0$ when $\ep\to0$. It is easy to
see that $D_3\le (\sum_{i>I} \be_i) (\|\vfi\|_0+\|\vfi\|_0)\le
2\de\|\vfi\|_0$.
If we take $\ep_0\in]0,\ep(I)[$ such that $|\be_i^\ep-\be_i|<\de/I$
for $\ep\in]0,\ep_0[$ and $i=1,\ldots,I$, then $\sum_{i=I+1}^{l(\ep)}
\be_i^\ep<2\de$ because $\sum_{i=1}^{l(\ep)}=1$.
Thus $D_2\le \sum_{i=I+1}^{l(\ep)}
(\be_i^\ep+\be_i)|\mu_i^\ep(\vfi)|\le
3\de\|\vfi\|_0$ for $\ep\in]0,\ep_0[$.
This shows that $D(\vfi,\ep)\to0$ when $\ep\to0$ for any fixed
continuous map $\vfi:M\to\real$ and concludes the proof of Theorem~\ref{thmA}.

\subsection{Convergence to finitely many measures}
\label{finitemeasures}
To prove Theorem~\ref{thmB} we now assume that $f$ satisfies
Conditions $\un B$ and $\un C$ for a finite family
$\La_1,\ldots,\La_N$, of attractors.
We observe that although Lemma~\ref{basin=stable} no longer holds, we
may still define the thresholds $\ep_1\ge\ldots\ge\ep_N>0$ as before
for a physical random perturbation $F$ of $f$ by
Proposition~\ref{uniqphysical}.
For $\ep\in]0,\ep_N[$ and $x\in M$ there are $l=l(\ep)\ge N$ physical
probability measures such that
$\mu^\ep(x)=\al_1^\ep(x)\mu_1^\ep+\ldots+\al_{l}^\ep(x)\mu_{l}^\ep$
where $\al_1^\ep(x),\ldots,\al_l^\ep(x)\ge0$ and $\sum_{i=1}^{l(\ep)}
\al_i^\ep(x)=1$, by Propositions~\ref{asymptotic} and~\ref{decomposition}.

If we show that $l(\ep)=N$ for $\ep\in]0,\ep_N[$, then Condition $\un
B$ and Proposition~\ref{uniqphysical} together imply that
$\mu_i^\ep\stackrel{w^*}{\longrightarrow}\mu_i$
when $\ep\to0$ for $i=1,\ldots, N$
and so we would get~\refer{thmB1}.

Arguing by contradiction, we suppose that $l(\ep)>N$ for
$\ep\in]0,\ep_N[$.  Then there is a physical
  measure $\mu_{N+1}^\ep$ and $\supp(\mu_{N+1})$ is a trapping region
  for $f^K$, where $K=K(\ep)$ from Condition II on $F$, by
  Theorem~\ref{thm0}. According to Proposition 2.9 of~\cite{Sh}, any
  trapping region for $f^K$ contains a trapping region for $f$. Then
  Condition $\un C$ would ensure that $\supp(\mu_{N+1}^\ep)$ contains
  $\La_i$ for $i\in\{1,\ldots,N\}$.
But we note that $\supp(\mu_i^\ep)\supset \La_i$ for
$i=1,\ldots,N$, by Proposition~\ref{uniqphysical} since $\ep<\ep_N$ and
also remark that $(\supp(\mu_i^\ep))_{i=1,\ldots,l(\ep)}$ are
  pairwise disjoint by Theorem~\ref{thm0}. We have arrived at a
  contradiction. Hence we have showed that $l(\ep)=N$ for all
  $\ep\in]0,\ep_N[$ and thus proved Theorem~\ref{thmB}.


{\em
V\'{\i}tor Ara\'ujo

Centro de Matem\'atica da Universidade do Porto

Rua do Campo Alegre  687, 4169-007 Porto, Portugal.

{\rm and}

Instituto de Matem\'atica,  Universidade Federal do Rio de Janeiro,

C. P. 68.530,  21.945-970 Rio de Janeiro, RJ-Brazil

Email: vitor.araujo@im.ufrj.br {\rm and} vdaraujo@fc.up.pt

}

\end{document}